\def\ZZ         {{\bf Z}}

\def\CC         {{\bf C}}
\def\QQ         {{\bf Q}}
\def\PP         {{\bf P}}
\def\ii         {{\rm i}}
\def\ee         {{\rm e}}

\def\dim        {{\rm dim}}
\def\sin        {{\rm sin}}

\def\Ell        {{\cal ELL}}

\documentstyle[twoside,12pt]{article}

\newtheorem{prop}{Proposition}[section]
\newtheorem{dfn}[prop]{Definition}   
\newtheorem{theo}[prop]{Theorem}

\newtheorem{rem}[prop]{Remark} 
 
\newtheorem{lem}[prop]{Lemma}

\begin{document}

\title{Higher elliptic genera} 

\author{Lev Borisov \\
\small Department of Mathematics\\
\small University of Wisconsin, Madison WI 53706\\
\bigskip
\small email: borisov@math.wisc.edu \\
         Anatoly Libgober \\
\small Department of Mathematics \\
\small University of Illinois, Chicago, IL 60607 \\
\small email: libgober@math.uic.edu}


\maketitle

\begin{abstract}
We show that elliptic classes introduced in \cite{annals} for 
spaces with infinite fundamental groups yield Novikov's type 
higher elliptic genera which are invariants of K-equivalence.
This include, as a special case, the
birational invariance of higher Todd classes 
studied recently by J.Rosenberg and 
J.Block-S.Weinberger.
We also prove the  modular properties of these genera, show that 
they satisfy a McKay correspondence, and consider their twist by 
discrete torsion.
\end{abstract}

\section{Introduction and statements of results.}

In works \cite{R} and \cite{BW} the authors considered the birational
properties of higher Todd, $L$ and $\hat A$ ``genera'' for nonsingular 
algebraic varieties. These ``genera'' are originated in the Novikov's 
conjecture that for a compact closed 
manifold $X$ with a fundamental 
group $\pi$ the higher signatures 
$\sigma_{\alpha}:( L(X) \cup f^*(\alpha))[X]$ are 
homotopy invariants. Here $L(X) \in H^*(X,\QQ)$ is the total $L$-class, 
$f$ is the map from $X$ to the classifying space $B\pi$ of the group 
$\pi$, $\alpha \in H^*(B\pi, \QQ)$ and $[X]$ is the fundamental class of 
$X$ (cf. \cite{JDavis} for a survey of Novikov's conjecture).    
The work \cite{R} conjectures that in the case when 
$X$ is projective algebraic and nonsingular the higher Todd invariants
$(Td(X) \cup f^*(\alpha))[X]$, where $Td$ is the total Todd class, 
are {\it birational} invariants. The work \cite{BW} contains 
a proof of this conjecture but also raises the problem of generalizing 
it to higher elliptic genera. 
The purpose of this note 
is to present such a generalization for the (two-variable) elliptic genus.
In fact, this generalization is done here in a wider category
consisting of Kawamata log-terminal pairs $(X,D)$ for which $X$ supports 
the action of a finite group $G$ so that the pair $(X,D)$ is $G$-normal
(cf. below). In particular we obtain the invariance of elliptic 
genus for $K$-equivalences, which, due to the fact that the Todd genus
is a specialization of the elliptic genus, contains the birational
invariance of higher Todd genus as a special case. 
We show also that many properties of 
ordinary elliptic genus (indicated in the abstract above)
also are satisfied by the higher 
elliptic genera. These two variable invariants extend
one variable higher elliptic genera considered (in smooth case) 
in connection with rigidity properties in \cite{liu}.

Let $X$ be a complex manifold. The (two variable) 
elliptic genus is defined as the holomorphic Euler characteristic
of

\begin{equation}
y^{-{{dim X} \over 2}} \otimes_{n \ge 1} (\Lambda_{-yq^{n-1}}
T^*_X \otimes \Lambda_{-y^{-1}q^n} T_X \otimes S_{q^n} T^*_M
\otimes S_{q^n}T_M)
\end{equation}
where $T_X$ is the tangent bundle, $T_X^*$ is its dual and for a bundle
$V$, the element $S_tV$ (resp. $\Lambda_tV$) is the power series 
over the semigroup of vector bundles on $X$  given by 
$\sum_{k \ge 0} t^kS^k(V)$ (resp. $\sum_{k \ge 0} t^k\Lambda^k(V)$).
One views this holomorphic Euler characteristic 
as a function on $\CC \times H$ where $H$ is the 
upper half plane using $y=exp(2 \pi i z), q=exp(2 \pi i \tau), 
z \in \CC, \tau \in H$. As such, it becomes a holomorphic function 
on $\CC \times H$. Its value at $q=0,y=0$ (i.e. the limit 
when $\Im z, \Im \tau \rightarrow \infty$), is the Todd genus of $X$,
as is seen directly from 
the definition.
It has the $L$ and $\hat A$ genera as certain limit values as well
(cf. \cite{invent}).

It follows from the Riemann-Roch theorem and a direct calculation 
(cf. \cite{invent}) that
the elliptic genus can be written in terms the Chern roots $x_i$ of the 
tangent bundle of $X$ as $\Ell(X)[X]$ where

\begin{equation}\label{ellipticclass}
\Ell(X)=\prod_i x_i {{\theta ({{x_i} \over {2 \pi \ii}}-z,\tau)} \over
{\theta ({{x_i} \over {2 \pi \ii }}, \tau)}} 
\end{equation}
and  
\begin{equation}
\theta(z,\tau)=q^{1 \over 8}  (2 \sin \pi z)
\prod_{l=1}^{l=\infty}(1-q^l)
 \prod_{l=1}^{l=\infty}(1-q^l \ee^{2 \pi \ii z})(1-q^l \ee^{-2 \pi \ii
z})
\end{equation}
is the classical theta function considered as the series in $y,q$
(cf. \cite{Chandra}).

Now let $X$ be a complex manifold as above and let $\pi$ be its 
fundamental group. Let  
$f: X \rightarrow B\pi$ be the corresponding map  
and let $\alpha \in H^k(\pi,\QQ)$. 

\begin{dfn}\label{higherelliptic}
The higher elliptic genus is 
$$Ell_{\alpha}(X)=(\Ell(X) \cup f^*(\alpha))[X]$$
where the elliptic class $\Ell(X) \in H^*(X,\QQ)$ is given by 
(\ref{ellipticclass}).
\end{dfn} 

Modular property of $Ell_{\alpha}(X)$ is described in 
theorem \ref{modularity} below but first we shall outline 
the extension to non simply-connected case of
the generalizations
of elliptic genus introduced in \cite{duke} and \cite{annals}.

Let $X$ be a normal projective algebraic variety and $D=\sum a_iD_i$ 
be a linear combination of distinct irreducible divisors with rational 
coefficients. The pair $(X,D)$ is called Kawamata log-terminal 
(cf. \cite{KMM}) if $K_X+D$ is $\QQ$-Cartier and there
is a birational morphism $f: Y \rightarrow X$ such that the union of 
the proper preimages of components of $D$ and the components of 
exceptional set  $E=\cup E_j$
 form a normal crossing divisor  
such that $K_Y=f^*(K_X+\sum a_iD_i)+\sum \alpha_jE_j$ 
and $\alpha_j >-1$
(here 
$K_X,K_Y$ are the canonical classes of $X$ and $Y$ respectively).
The triple $(X,D,G)$ where $X$ is a nonsingular variety, 
$D$ is a divisor and $G$ is a finite group of biholomorphic automorphisms is  
 called $G$-normal (cf. \cite{bat}, \cite{annals}) 
if the components of $D$ form a normal crossings divisor 
and the isotropy group of any point acts trivially on the 
components of $D$ containing this point.

\begin{dfn}(cf. \cite{annals} definition 3.2)\label{orbifoldelliptic}
 Let $(X,D)$ be a Kawamata log terminal $G$-normal pair
and $D=-\sum \delta_kD_k$. The {\it orbifold} elliptic class of 
 $(X,D,G)$ is the class in $H_*(X,\QQ)$ given by:

$$
{\Ell}(X,D,G;z,\tau):=
\frac 1{\vert G\vert }\sum_{g,h,gh=hg}\sum_{X^{g,h}}[X^{g,h}]
\Bigl(
\prod_{\lambda(g)=\lambda(h)=0} x_{\lambda} 
\Bigr)
$$
$$\times\prod_{\lambda} \frac{ \theta(\frac{x_{\lambda}}{2 \pi \ii }+
 \lambda (g)-\tau \lambda(h)-z )} 
{ \theta(\frac{x_{\lambda}}{2 \pi \ii }+
 \lambda (g)-\tau \lambda(h))}  \ee^{2 \pi \ii \lambda(h)z}
$$
$$\times\prod_{k}
\frac
{\theta(\frac {e_k}{2\pi\ii}+\epsilon_k(g)-\epsilon_k(h)\tau-(\delta_k+1)z)}
{\theta(\frac {e_k}{2\pi\ii}+\epsilon_k(g)-\epsilon_k(h)\tau-z)}
{}
\frac{\theta(-z)}{\theta(-(\delta_k+1)z)} \ee^{2\pi\ii\delta_k\epsilon_k(h)z}.
$$
\end{dfn}
where $X^{g,h}$ denotes an irreducible component of the fixed set of the 
commuting elements $g$ and $h$ and $[X^{g,h}]$ denotes the image 
of the fundamental class in $H_*(X)$. The restriction of $TX$ to $X^{g,h}$
splits into linearized bundles according to the ($[0,1)$-valued) characters 
$\lambda$  of $\langle g,h\rangle$. 
Moreover, $e_k=c_1(E_k)$ 
and $\epsilon_k$ is the character of ${\cal O}(E_k)$ restricted to $X^{g,h}$ 
if $E_k$ contains $X^{g,h}$ and is zero otherwise.

\begin{dfn}\label{higherclasstriple}
 Let $\pi$ and $\alpha$ be as in definition \ref{higherelliptic}
and other notations as in the definition \ref{orbifoldelliptic}.
Then 
$$Ell_{\alpha}(X,D,G)=(\Ell(X,D,G) \cap f^*(\alpha))_0$$
i.e. the component of degree zero of the class in $H_*(X,\QQ)$.
\end{dfn}

Special cases of the elliptic class $\Ell(X,D,G)$ 
are the following. 

a) If $G$ is the trivial group one obtains higher elliptic genus  
of Kawamata log terminal pairs $\Ell(X,D)$. If $D=\sum \delta_kD_k$ 
and $d_k \in H^2(X,\QQ)$ is the cohomology class dual to $D_k$ then:

\begin{equation}\label{ellipticclasspair}
\Ell(X,D)=(\prod_l 
\frac{(\frac{x_l}{2\pi\ii})\theta(\frac{x_l}{2\pi\ii}-z)\theta'(0)} 
{\theta(-z)\theta(\frac{x_l}{2\pi\ii})} 
\Bigr)\times 
\Bigl(\prod_k 
\frac{\theta(\frac{d_k}{2\pi\ii}-(\delta_k+1)z)\theta(-z)} 
{\theta(\frac{d_k}{2\pi\ii}-z)\theta(-(\delta_k+1)z)} 
\Bigr) 
\end{equation}

b) If $D=\emptyset$ one obtains the higher orbifold elliptic genus
$Ell_{orb,\alpha}$ which is the value on $[X]$ 
of the cup product with $f^*(\alpha)$ of the orbifold 
elliptic class:

\begin{equation}\label{orbifoldclass}
 \Ell_{orb}(X,G)=
{1 \over {\vert G \vert }}
\sum_{g,h, gh=hg} 
\Bigl(
\prod_{\lambda(g)=\lambda(h)=0} x_{\lambda} 
\Bigr)
\prod_{\lambda} {{ \theta(\tau,{{x_{\lambda}} \over {2 \pi \ii }}+
 \lambda (g)-\tau \lambda(h)-z )} \over 
{ \theta(\tau,{{x_{\lambda}} \over {2 \pi \ii }}+
 \lambda (g)-\tau \lambda(h))} } e^{2 \pi \ii \lambda(h)z}
[X^{g,h}]
\end{equation}

Subscript $\alpha$ will denote the twisting by the class $\alpha$ 
of the genus corresponding to 
(\ref{ellipticclasspair}) and (\ref{orbifoldclass}). 
Either of these is a special case of \ref{higherclasstriple}.

In the case $D=\emptyset, G=\{1\}$ we obtain the class given by  
(\ref{ellipticclass}) and the higher elliptic genus defined 
in \ref{higherelliptic}.

Finally recall the following (cf. \cite{zagier})
\begin{dfn}
A Jacobi form of index $t \in {1 \over 2}\ZZ$ and weight $k$ is a holomorphic 
function $\chi$ on $H \times \CC$ satisfying the following functional 
equations:
$$
\chi({{a\tau+b} \over {c\tau+d}},{z \over {c\tau+d}})=
(c\tau+d)^ke^{{2 \pi i t c z^2} \over {c\tau+d}}\chi(\tau,z)$$ 
$$\chi(\tau,z+\lambda\tau+\mu)=(-1)^{2t(\lambda+\mu)}
e^{-2\pi i t(\lambda^2\tau+2 \lambda z)}\chi(\tau,z)
$$
\end{dfn}

Important property of the elliptic genus is that for a 
$SU$-manifold the elliptic genus 
is a Jacobi form having weight zero and index ${{\rm dim} X \over 2}$.

\begin{theo}\label{modularity}
 Let $X$ be a $SU$-manifold,  $d=\dim X$,
$\pi=\pi_1(X)$ and $\alpha \in H^{2k}(\pi,\QQ)$. Then the higher elliptic genus 
$(\Ell(X) \cup f^*(\alpha))[X]$ is a Jacobi form having index $d \over 2$ 
and weight $-k$. It has the Novikov signature, the higher 
Todd genus and higher $\hat A$-genus as specializations. 

More generally, let $X,D$ be a Kawamata log-terminal $G$-normal pair
where $G$ is a finite group. If $m(K_X+D)$ is a trivial Cartier 
divisor,  $n$ is the order of the image  
$G \rightarrow Aut H^0(X,m(K_X+D))$ and $\alpha \in H^k(\pi_1(X),\QQ)$ 
as above then 
$\Ell_{\alpha}(X,D,G,z,\tau)$ is a Jacobi form 
having weight $-k$ and the index ${{\dim X}\over 2}$ for the 
subgroup of Jacobi group generated by:
$$(z,\tau) \rightarrow (z+mn,\tau),(z,\tau) \rightarrow (z+mn \tau),
(z,\tau) \rightarrow (z,\tau), (z,\tau) \rightarrow ({z \over \tau},
-{ 1\over {\tau}})$$
\end{theo}

Next recall that two manifolds $X_1,X_2$ are called $K$-equivalent if 
there is a smooth manifold $\tilde X$ and a diagram:

\begin{equation}\label{kequivalencepicture}
 \matrix{ &  & \tilde X & & & \cr
                         & \phi_1\swarrow &  & \searrow \phi_2 & & \cr
                         X_1 & & & & X_2 \cr}   
\end{equation}
in which $\phi_1$ and $\phi_2$ are birational  
morphisms and $\phi_1^*(K_{X_1})$ and $\phi_2^*(K_{X_2})$ 
are linearly equivalent.

\begin{theo}\label{kequivalence} For any $\alpha \in H^*(B\pi,\QQ)$ 
the higher elliptic genus $(\Ell(X) \cup f^*(\alpha),[X])$
is an invariant of $K$-equivalence. Moreover, if $(X,D,G)$ 
and  $(\hat X, \hat D, G)$ 
are $G$-normal and Kawamata log-terminal and if 
$\phi: (\hat X,\hat D) \rightarrow (X,D)$ is $G$-equivariant
such that 
$$\phi^*(K_X+D)=K_{\hat X}+\hat D$$ then 
$$Ell_{\alpha}(\hat X,\hat D,G)=Ell_{\alpha}(X,D,G)$$
In particular the higher elliptic genera (and hence the higher
signatures and $\hat A$-genus)
are invariant for  crepant 
birational morphisms and 
the specialization into the Todd class
is birationally invariant.
\end{theo}

\begin{rem} Since 
the fundamental groups are rather restricted by algebraic geometry,
one may wonder when higher genera yield new invariants. 
This also depends on the existence
of nontrivial cohomology classes of the fundamental group
(cf. remark \ref{smallcodclasses} below). 
By a result of Beauville 
a fundamental group of a Calabi Yau manifold is an extension of 
a free abelian group by a finite group (cf. \cite{beauville}) 
so one does obtain new invariants if the rank of this abelian 
group is positive. On the other hand, the higher elliptic genera 
of pairs are defined with $X$ being arbitrary projective algebraic 
manifold and class of groups from which one obtains higher invariants
is much bigger than in Calabi Yau case.
\end{rem}

\section{Proofs}

{\it Proof of the theorem \ref{modularity}}. The argument is 
essentially the same as for the special ($k=0$) case 
in \cite{invent} or in more general case of orbifold elliptic genus of 
pairs dealt with in \cite{annals}. We spell out the argument 
only in the first case of the theorem \ref{modularity}. It is enough to 
check the transformation formulas on the generators
of the Jacobi groups i.e. that 

\begin{equation} 
\chi(M,z,\tau+1)=\chi(M,z,\tau) 
\label{modular1} 
\end{equation}
\begin{equation}
\chi(M,{z \over \tau},-{1 \over \tau})=\tau^k \ee^{{ \pi \ii d z^2} \over \tau}
\chi(M,z,\tau) \label{modular2} 
\end{equation}
\begin{equation} \chi(M,z+\tau,\tau)= (-1)^d
\ee^{- \pi \ii d (\tau + 2z)} \chi(M,z,\tau) \label{modular3}
\end{equation}
\begin{equation}
\chi(M,z+1,\tau)=(-1)^d \chi(M,z,\tau) .\label{modular4} \end{equation}

Let 
\begin{equation}
\prod x_i  {{\theta({x_i \over {2 \pi \ii }}-z, \tau)}
\over {\theta ({{x_i} \over {2 \pi \ii }},\tau)}}= 
\sum_{\bf k} Q_{\bf k}(z,\tau) {\bf x}^{\bf k}
\end{equation}
where $\bf x$ is a product of powers of $x_i$ and $\bf k$ is
multiindex. The higher elliptic genus is the linear 
combination of functions $Q_{\bf k}(z,\tau)$ with 
coefficients $f^*(\alpha) \cup {\bf x}^{\bf k}[X]$ for all 
${\bf k}, \ \vert {\bf k} \vert =d-k$.
The transformation formulas for $\theta$-function yield that the left hand side
satisfies (\ref{modular1}),(\ref{modular2}) and (\ref{modular3}).
Hence these relations are valid for all $Q_{\bf k}(z,\tau)$.
For the transformations $\tau \rightarrow {-{1 \over {\tau}}},
x_i \rightarrow {{x_i} \over \tau}, z \rightarrow {z \over \tau}$
we have:

$$
\sum_{\bf k} Q_{\bf k}({z \over \tau},-{1 \over {\tau}})
 ({{\bf x} \over \tau})^{\bf k}=
\prod_i ({x_i \over {\tau}}) {{\theta (-{z \over
\tau}+{{x_i}
\over {2 \pi \ii \tau}},-{1 \over \tau})} \over
{\theta ({{x_i} \over {2 \pi \ii \tau}},-{1 \over \tau})}}=
({1 \over \tau})^d \prod \ee^{-{z x_i} \over \tau}x_i
{{\ee^{{\pi \ii z^2} \over \tau}
\theta(-z+
{{x_i} \over {2 \pi \ii }},\tau )} 
\over {\theta ({{x_i} \over {2 \pi \ii }},\tau )}}$$
\begin{equation}
=({1 \over \tau})^d \prod x_i{{\ee^{{\pi \ii z^2} \over \tau}
\theta(-z+
{{x_i} \over {2 \pi \ii }},\tau )} 
\over {\theta ({{x_i} \over {2 \pi \ii }},\tau )}}=
\sum_{\bf k}({1 \over \tau})^d {\ee^{{\pi \ii dz^2} \over \tau}}
Q_{\bf k}(z,\tau) {\bf x}^{\bf k}
\end{equation}
(the equality before the last is due to the vanishing of $c_1$).
Hence the coefficient of a monomial ${\bf x}^{\bf k}$ with 
$\vert {\bf k} \vert =d-k$ satisfies:
$$Q_{\bf k}({z \over \tau},-{1 \over {\tau}})=
t^{-k}{\ee^{{\pi \ii dz^2} \over \tau}}Q_{\bf k}(z,\tau)$$

\begin{rem}\label{smallcodclasses} Elliptic genus $Ell_{\alpha}(X)=0$ if 
$\alpha \in H^{2k}(\pi_1,\QQ), k>d$
and is a multiple of the Jacobi form:
$$({\theta(\tau,z) \over {\theta'(0,\tau)}})^d$$   
for $\alpha \in H^{2d}(\pi_1,\QQ)$.
The latter has weight $-d$ and index $d \over 2$. 

For a class 
$\alpha \in H^{2d-2}(\pi_1(X),\QQ)$ the genus $Ell_{\alpha}(X)$
is determined by the higher Todd genus but elliptic genera 
corresponding to classes of large codimension cannot  
be described in terms of 
higher Todd (or $\chi_y$)-genus (cf. \cite{invent}, theorem 2.7). 
\end{rem}

\bigskip
\noindent {\it Proof of the theorem \ref{kequivalence}}. We
apply Theorem 3.5
of \cite{annals} to 
the resolution $\phi_1:\tilde X\rightarrow X_1$ 
of the diagram (\ref{kequivalencepicture}) to get
a direct image formula
\begin{equation}
(\phi_1)_*{\Ell}(\tilde X,\tilde D,G;z,\tau)={\Ell}(X_1,D_1, G;z,\tau).
\end{equation}
Here $\tilde D$ is defined as usual as having the same coefficients 
as $D_1$ at components of the proper preimage of $D_1$ and having 
the coefficients at the exceptional divisors of $\phi_1$ determined
from $\phi^*(K_X+D)=K_{\hat X}+\tilde D$.

Together with the diagram
\begin{equation}
 \matrix{~~\tilde X & & \cr
        \phi_1 \hskip -.2cm
 \ \downarrow & \buildrel \tilde f \over 
\searrow & \cr 
            ~~X_1 & \buildrel f_1 \over \rightarrow  & B\pi}
\end{equation}
the direct image formula yields
$$ (\Ell_{\alpha}(\tilde X,\tilde D, G) \cap \tilde f^*(\alpha))_0=$$
\begin{equation}
(\phi_1)_*(\Ell(\tilde X, \tilde D,G) \cap \hat (\phi_1 \circ 
f_1)^*(\alpha))_0=
(\Ell(X_1, D_1,G) \cap f_1^*(\alpha))_0.
\end{equation} 
The $K$-equivalence means that the divisor $\tilde D$ calculated
for $\phi_1$ is the same as the divisor $\tilde D$ calculated for 
$\phi_2$, which shows that higher elliptic genera 
are  invariant under $K$-equivalences.
The claim about crepant morphisms 
is then immediate. Finally, 
as one sees from formula (\ref{ellipticclasspair}), 
in the limit $\Im z \rightarrow \infty, \tau \rightarrow \infty$ the 
elliptic class of pair loses its dependence on $D$. Hence the push forward
formula is valid without assumption on the canonical class and  
the higher Todd classes are invariant for arbitrary birational 
morphisms and not just $K$-equivalences.

\section{Further properties of higher genera}

\subsection{Higher elliptic genera of singular varieties.}

One of the consequence of previous discussion is existence of 
well defined higher elliptic genera of projective algebraic 
varieties with log-terminal singularities.


If $\tilde X$ is a resolution of a projective variety $X$ then 
$\pi_1(X) \ne \pi_1(\tilde X)$ 
in general (for example image of generic projection of $X$ in 
$\PP^{\dim X+1}$ is simply-connected for any $X$ (cf. \cite{Fulton}).
However we have the following result due to Takayama (cf. \cite{takayama}(*)):
\footnote{(*)We thank C.Hacon and J.McKernan for pointing out 
the reference and J.McKernan for further comments.} 

\begin{lem}\label{takayamaresult}
 Assume that $X$ has only log-terminal singularities 
(or more generally $(X,\Delta)$ has divisorial Kawamata log-terminal 
singularities in terminology of \cite{hacon}) and let $f: X' \rightarrow X$
be a resolution of singularities of $X$. Then $\pi_1(X')=\pi_1(X)$.
\end{lem}

It follows from the theorem \ref{kequivalence} that 
the following definition yields result independent of 
resolution.

\begin{dfn}Let $X$ be a projective algebraic 
variety with $\QQ$-Gorenstein log-terminal singularities.
Let $\alpha \in H^*(\pi_1(X),\QQ)$ be the cohomology 
class of its fundamental group. If $\phi: \tilde X \rightarrow X$ 
is a resolution 
of singularities of $X$, $K_{\tilde X}=\phi^*(K_X)+\tilde D$ and 
$\alpha$ is viewed as the 
element in the $H^*(\pi_1(\tilde X),\QQ)\buildrel {\rm dfn} 
\over = H^*(\pi_1(X),\QQ)$
then:
$$Ell_{\alpha}(X)=Ell_{\alpha}(\tilde X, \tilde D)$$
\end{dfn}

Note that similarly one can define the $Ell_{\alpha}(X,D)$ 
where $(X,D)$ is such that $K_X+D$ is $\QQ$-Cartier and
such that $(X,D)$ is a log-terminal pair. 
A large class of varieties with singularities as in this definition 
can be obtained by looking at
 the quotients of nonsingular varieties by the action of a finite
group acting via biholomorphic transformations. In the next section
we show that the calculation of the higher 
elliptic genus of quotients 
can be done in terms the action on nonsingular variety. 
This extends the results of \cite{annals} 
for ordinary elliptic genus which 
in turn extend the results on Euler characteristics, 
$\chi_y$-characteristic etc. of quotients (see this reference
for review of the preceding work).

\subsection{Higher McKay correspondence}

Let $X$ be nonsingular and $G$, as before, a finite group 
of biholomorphic transformations. Let $\pi=\pi_1(X/G)$
and let $\alpha_{X/G} \in H^*(\pi_1(X/G),\QQ)$. Let 
$\mu_G: \pi_1(X) \rightarrow \pi_1(X/G)$ be the homomorphism 
corresponding to $X \rightarrow X/G$, 
$\alpha_X=\mu_G^*(\alpha_{X/G}) \in H^*(\pi_1(X),\QQ)$.  
The next theorem describes the invariant $Ell_{orb,\alpha}(X,G)$ 
given by class (\ref{orbifoldclass})
in terms of resolution of singularities of $X/G$ 
(which is the classical McKay correspondence 
between the Euler characteristic of minimal resolution and 
the order of the group in the case $\pi_1=\{1\},q=0,
y=1,X=\CC^2,G \subset SL_2(\CC)$). In the case $\alpha=1$ 
the McKay correspondence for elliptic genera was 
conjectured in \cite{duke} and proven in \cite{annals}.
In the case of arbitrary $\alpha$ we have the following.

\begin{theo} Let $X$ be a nonsingular projective variety, 
$G$ acts biholomorphically on $X$, $\mu: X \rightarrow X/G$,
$D=\sum (\nu_i-1)D_i$ is the ramification divisor of $\mu$ 
and $\Delta_{X/G}=\sum {{\nu_i-1} \over {\nu_i}}\mu(D_i)$.
Then:
$$Ell_{orb,\mu^*(\alpha_X/G)}(X,G;z,\tau)=
{Ell}_{\alpha_{X/G}}(X/G,\Delta_{X/G};z,\tau)$$
\end{theo}  
Indeed, this follows by the same argument as 
the one used in the proof of the theorem 5.3 in \cite{annals}
by applying (obtained in lemma 5.4, \cite{annals})
the push forward formulas and the 
projection formula to the diagram:

\begin{equation} \matrix{\mu_Z: & \hat Z & \rightarrow & \hat X \cr
                               & \phi \downarrow & & \downarrow    \cr
                         \mu:   &   X   & \rightarrow & X/G  \cr}     
\end{equation} 
since both sides of the claimed equality are pushforwards 
of the same class $\Ell(\hat Z, \hat D,G)$ on $\hat Z$.

\subsection{Concluding remarks: flops, rigidity and discrete torsion}

Higher elliptic genera yield 
a description of $\Omega^{SU}(B\pi) \otimes \QQ $ modulo flops
using Jacobi forms 
extending Totaro's (cf.\cite{totaro})
description in the case $\pi=1$.

\begin{dfn}
We shall say that two maps $f_X: X \rightarrow B\pi$ and 
$f_{X'}: X' \rightarrow B\pi$ are differ by a flop if: 

a) there is a map of an (almost) 
complex  space $f_Y: Y \rightarrow B\pi$ such that the singular set 
${\rm Sing} Y= Z$  is a manifold having in $Y$ 
a regular neighborhood  biholomorphic to 
$Z \tilde \times \cal V$ where ${\cal V} \subset {\CC}^4$ is given by 
$xy=zw$ and

b) there are maps $\pi_X: X \rightarrow Y$ and $\pi_{X'}: X' \rightarrow Y$
such that $\pi_X$ and $\pi_X'$ are small resolutions of 
$Y$ (i.e. $\pi_{X}^{-1}(Z) \rightarrow Z$  and 
$\pi_{X'}^{-1}(Z) \rightarrow Z$
are $\PP^1$ fibrations yielding at each point of $Z$ two 
different $\PP^1$-resolutions of the node $\cal V$)
and such that one has the commutative diagram:
$$\matrix{ &  & Y & & \cr
           & \pi_X \nearrow & & \pi_{X'} \nwarrow  & \cr
         X & & f_Y\downarrow & & X' \cr 
             & f_X \searrow & & f_{X'} \swarrow  & \cr 
         &  & B\pi & & \cr  }$$ 

\end{dfn}

\begin{prop}
Let $I_{\pi}$ (resp. $I$) be the ideal in $\Omega^{U}(B\pi)$ generated by the 
differences $(X,\pi_X)$ and $(X',\pi_{X'})$ from the above definition 
(resp. $X'-X$ where $X'$ and $X$ 
are differ by a classical flop (cf. \cite{totaro}).)
Then one has the following isomorphism of 
graded $\Omega_*^U$-modules:
$$\Omega_*^U(B\pi)\otimes \QQ/I_{\pi}=H_*(B\pi,\QQ) \otimes_{\QQ} \Omega_*^U/I=
H_*(B\pi,\QQ) \otimes \QQ[x_1,x_2,x_3,x_4]$$
(cf. \cite{totaro} for geometric description of the isomorphism:
$\Omega_*^U/I=\QQ[x_1,x_2,x_3,x_4]$).
In particular:
\begin{equation}\label{cobordismiso}
Hom (\Omega^{SU}_d(B\pi)/I_{\pi} \cap \Omega^{SU}_d(B\pi),\QQ) =
\oplus_{k \in 2\ZZ} H^k(B\pi,{\rm Jac}_{-k,{{d} \over2}}) 
\end{equation}
where ${\rm Jac}_{-k,{{d} \over2}}$ is the space of Jacobi forms
having weight $-k$ and index $d \over 2$ (*).
\footnote{(*)The isomorphism  (\ref{cobordismiso}) is a counterpart of 
the theorem 4.3 in \cite{R}. The assumption in \cite{R}
that $H_k(B\pi,\QQ)$ 
is spanned by classes of maps form projective varieties can be 
omitted if birational invariance is stated as triviality on 
the almost complex manifolds which are projectivisations 
of complex bundles: indeed any birational equivalence is 
composition of blowups and blowdowns (cf. \cite{akmw}) and 
the difference of a manifold and its blow in $\Omega^U$ is 
a projectivised bundle (cf. \cite{hitchin}).} 
\end{prop}
This follows from the isomorphism (cf.\cite{conner}):
\begin{equation}\label{moduloflops}  
\oplus H_k(B\pi,\QQ) \otimes \Omega^{U}_{d-k/2}\otimes \QQ 
\rightarrow \Omega^{U}_{d}(B\pi)\otimes \QQ
\end{equation} 
which assigns to an almost complex manifold $\pi_N: N \rightarrow B\pi$ 
representing a homology class $\alpha$ and an almost complex manifold
$M$ the map of $M \times N$ which is the composition of the projection on $N$ 
and $\pi_N$. The condition that $X'-X \in I$ 
is equivalent to 
$X'-X=
F \tilde \times Z \rightarrow B\pi$ 
where $F$ is the almost complex $SU$-manifold
which is homological $\CC\PP^3$ (cf. \cite{totaro}). The map
(\ref{moduloflops}) takes it to zero. Converse it similar. 
Note that for $\pi_1=(1)$ (\ref{cobordismiso}) becomes the 
identification in \cite{totaro}.

Let $X$ be a Calabi Yau algebraic manifold with the fundamental group
$\pi$ and let $X \rightarrow Alb(X)$ be the Albanese map. Since the
fundamental group of $X$ has a free abelian group of rank 
$2 {\rm dim} Alb(X)$ as a subgroup of finite index 
(cf. \cite{beauville}) we have the following description of the
higher ellptic genera which can appear in algebraic case.
Using the isomorphism 
\begin{equation}\label{jacobi}
J_{-k,{d \over 2}}=\oplus J_{0,{{d-k} \over 2}}
\end{equation}
(which follows for example since 
the ring of Jacobi form has a system of generators such that 
only one generator  $({\theta(\tau,z) \over {\theta'(0,\tau)}})$ 
has negative weight (cf. \cite{gritsenko})) the isomorphism 
(\ref{moduloflops}) implies that the 
 higher genera coincide with the ordinary elliptic 
genera of preimages of (almost complex) 
submanifolds of the abelian variety $Alb(X)$ interpreted as 
elements of $J_{-k,{d \over 2}}$ using (\ref{jacobi}).

Other application of higher elliptic genera include:

\bigskip

\noindent a) Rigidity 
theorems extending Browder-Hsiang (cf. \cite{BH}) on 
rigidity of higher $\hat A$-genus on one hand and Bott-Taubes
results on rigidity of elliptic.
Case of higher one variable elliptic genus is discussed in 
\cite{liu} and \cite{gong} (cf. also \cite{waelder}).

\bigskip

\noindent
b) One has a version of higher elliptic genera twisted by 
discrete torsion extending the one considered in \cite{LS}. 
Let $X,D$ be a $G$-normal pair where $G$ is a finite 
group of biholomorphic automorphisms. A class $\nu \in H^2(G,U(1))$
defines $\delta(g,h) = \frac{\alpha(g,h)}{\alpha(h,g)}$ which allows 
to twist using $\nu$ the elliptic class in definition \ref{orbifoldelliptic}
to obtains the class: 

$$
{\Ell}^{\nu}(X,D,G;z,\tau):=
\frac 1{\vert G\vert }\sum_{g,h,gh=hg}\sum_{X^{g,h}}[X^{g,h}]\delta(g,h)
\Bigl(
\prod_{\lambda(g)=\lambda(h)=0} x_{\lambda} 
\Bigr)
$$
$$\times\prod_{\lambda} \frac{ \theta(\frac{x_{\lambda}}{2 \pi \ii }+
 \lambda (g)-\tau \lambda(h)-z )} 
{ \theta(\frac{x_{\lambda}}{2 \pi \ii }+
 \lambda (g)-\tau \lambda(h))}  \ee^{2 \pi \ii \lambda(h)z}
$$
$$\times\prod_{k}
\frac
{\theta(\frac {e_k}{2\pi\ii}+\epsilon_k(g)-\epsilon_k(h)\tau-(\delta_k+1)z)}
{\theta(\frac {e_k}{2\pi\ii}+\epsilon_k(g)-\epsilon_k(h)\tau-z)}
{}
\frac{\theta(-z)}{\theta(-(\delta_k+1)z)} \ee^{2\pi\ii\delta_k\epsilon_k(h)z}.
$$

This yields the following version of the elliptic genus:
\begin{equation}
Ell_{\alpha}^{\nu}(X,D,G)=(\Ell^{\nu}(X,D,G) \cap f^*(\alpha))_0
\end{equation}
This is a Jacobi form of weight $-k$ and index $d \over 2$. Such elliptic 
genus is also rigid for $S^1$ actions commuting with the action of $G$ and 
preserving $D$.

\bigskip

\bigskip               
\par  


\noindent 

\begin{thebibliography}{99}

\bibitem{akmw}  
D. Abramovich, K. Karu, K. Matsuki, J. W\l odarczyk,
{\em Torification and Factorization of Birational Maps}, 
J. Amer. Math. Soc. {\bf 15} (2002), no. 3, 531--572.

\bibitem{beauville} A.Beauville, {\em Vari\'et\'e Kahleriennes dont la
premiere class Chern est nulle}, J.Diff. Geom, (18), 1983, no.4 755-782. 


\bibitem{BW} J.Block, S.Weinberger, 
{\em Higher Todd classes and holomorphic group actions}.
preprint. math.AG/0511305

\bibitem{bat}V.Batyrev,
 Non-Archimedean integrals and stringy Euler numbers 
of log-terminal pairs.  J. Eur. Math. Soc. (JEMS)  1  (1999),  
no. 1, 5--33. 

\bibitem{invent} 
L.A. Borisov, A. Libgober, 
{\em Elliptic Genera of Toric 
Varieties and Applications to Mirror Symmetry}, 
Invent. Math. {\bf 140} (2000), no. 2, 453-485.


\bibitem{duke} L.A. Borisov, A. Libgober, 
{\em Elliptic genera of singular varieties},
Duke Math. J. {\bf 116} (2003),  no. 2, 319--351. 


\bibitem{annals} L.Borisov, A.Libgober, 
{\em McKay correspondence for 
elliptic genera},  Ann. of Math. (2)  161  (2005),  no. 3, 1521--1569.


\bibitem{BT} R.Bott, C.Taubes,
{\em On the rigidity theorems of Witten}.  
J. Amer. Math. Soc.  2  (1989),  no. 1, 137--186.  


\bibitem{BSY} J.P.Brasselet, J.Schuermann, S.Yokura, 
Hirzebruch classes and motivic
Chern classes for singular spaces, math.AG/0503492.



\bibitem{BH} W.Browder and W.C. Hsiang, 
{\em G-actions and the fundamental group},
Inventiones Math. 65, 1981/82 p.411-424. 

\bibitem{Chandra} K. Chandrasekharan, {\em Elliptic functions}, Fundamental
Principles of Mathematical Sciences, 
{\bf 281}, Springer-Verlag, Berlin-New York,
1985.


\bibitem{conner} P.Conner and E.Floyd, {\em Differentiable Periodic Maps}, 
Ergebnisse der Mathematic, band 33, Springer Verlag, 1964


\bibitem{JDavis} J.Davis, {\em Manifold aspects of the Novikov conjecture}, 
Surveys in surgery theory, vol.1, 195-224, Ann. of Math. Studies, 
Princeton Univ. Press. Princeton, N.J. 2000.

\bibitem{DLM} C. Dong, K. Liu, X. Ma, 
{\em On orbifold elliptic genus}, preprint math.DG/0109005.

\bibitem{gong} D.Gong, K.Liu, Rigidity of higher elliptic genera,
Ann. Global Anal.Geom. 14 (1996), 219-236.



\bibitem{zagier}
M. Eichler, D. Zagier, {\em The theory of Jacobi forms},
Progress in Mathematics, {\bf 55}, Birkhäuser Boston, Inc., Boston, Mass.,
1985

\bibitem{Fulton} W.Fulton, R.Lazarsfeld, {\em Connectctivity and its
applications in algebraic geometry}, Algebraic Geometry, (Chicago, Ill, 1980),
pp.26-92. Lecture Notes in Math. Springer Verlag, 1981.

\bibitem{GM} M.Goresky, R.MacPherson, {\em Stratified Morse Theory},
Springer Verlag, 1988.


\bibitem{gritsenko} V.Gritsenko, {\em Complex vector bundles and Jacobi forms},
Sriukaisekikenkyusho Kokyuroku No.1 1103 (1999), 71-85.


\bibitem{hacon} C.Hacon and J.McKernan, {\em Shokurov's Rational connectedness
conjecture}. math.AG/0504330

\bibitem{hitchin} N.Hitchin, {\em Harmonic Spinors}, 
Advances in Mathematics, (14) 1974, 1-55.


\bibitem{kollarm} J.Kollar,Y.Miyaoka, S.Mori, 
{\em Rationally connected varieties}, 
J.Algebraic Geometry, 1 (1992) no.3 p.429-448.

\bibitem{LS} A.Libgober, M.Szczesny,
{\em Discrete torsion, orbifold elliptic genera, 
and the chiral de Rham complex}, preprint, math.AG/0412422.


\bibitem{liu} Kefeng Liu, {\em On Mod 2 and Higher Elliptic Genera}, 
Comm. Math. Phys, 149, 71-95.

\bibitem{KMM} Y.Kawamata, K.Matsuda, K.Matsuki,
Introduction to the minimal model problem. 
 Algebraic geometry, Sendai, 1985,  
283--360, Adv. Stud. Pure Math., 10, North-Holland, Amsterdam, 1987. 


\bibitem{R} J.Rosenberg, {\em An analog of Novikov's conjecture in complex
algebraic geometry,preprint}, math.AG/0509526.
 

\bibitem{takayama} S.Takayama, {\em Local simple connectedness of resolutions
of log-terminal singularities}, Intern. Journal of Mathematics, vol.14 no 8.
p.825-836.

\bibitem{totaro}B. Totaro, {\em Chern numbers of singular varieties and
elliptic homology}, Ann. of Math. {\bf 151} (2000), no. 2, 757-791.

\bibitem{waelder} R.Waelder, Equivariant elliptic genera, math.AG/0603521.


\end{thebibliography}
\end{document}